\documentclass[a4paper]{amsart}

\usepackage{amssymb}
\usepackage{amsmath}
\usepackage{amsthm}
\usepackage{amsfonts}

\usepackage{a4wide}
\usepackage{graphicx}
\usepackage[usenames,dvipsnames]{xcolor}

\begin{document}

\includegraphics[width=6.0cm]{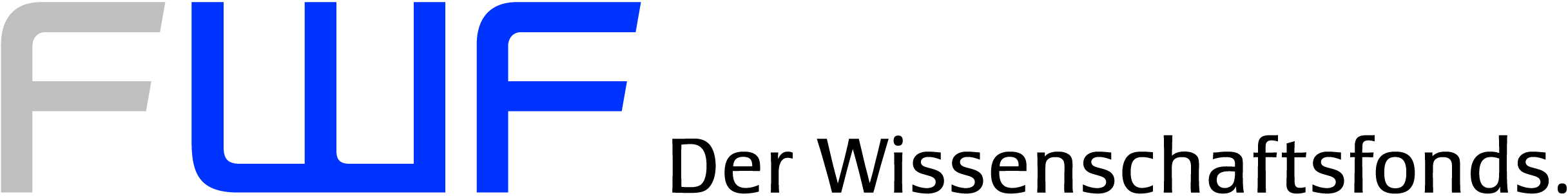}

\vspace{0.1in}

{\Huge

{\color{blue}
\begin{center}
\textbf{Group Theory of \\ Cyclic Cubic Number Fields}
\end{center}
}

\vspace{0.1in}

\renewcommand{\arraystretch}{1.4}

\begin{table}[ht]
\label{tbl:HeaderGTCCNF}
\begin{center}
\begin{tabular}{|ll|}
\hline
 \textbf{Conference:}             & 5th International Conference on \\
                                  & Mathematics and Computer Science \\
                                  & (MACOS) Bra\c sov 2024 \\
 \textbf{Place:}                  & Universitatea Transilvania din Bra\c sov \\
 \textbf{Venue:}                  & Bra\c sov (Kronstadt), Romania \\
 \textbf{Date:}                   & June 13 -- 15, 2024 \\
 \textbf{Author:}                 & Daniel C. Mayer (Graz, Austria) \\
 \textbf{Affiliation:}            & Austrian Science Fund (F{\color{blue}WF}) \\
\hline
\end{tabular}
\end{center}
\end{table}

\vspace{0.1in}

\begin{center}
A presentation within the frame of the \\
international scientific research project
\end{center}

\bigskip

{\color{blue}
\begin{center}
\textbf{Towers of \(p\)-Class Fields \\ over Algebraic Number Fields}
\end{center}
}

\bigskip

\begin{center}
supported by the Austrian Science Fund (F{\color{blue}WF}): \\
Projects J 0497-PHY and P 26008-N25, \\
and by the Research Executive Agency of the European Union (EUREA):\\
Project Horizon Europe 2021 -- 2027.
\end{center}

}

\newpage

{\normalsize

\noindent
Daniel C. Mayer (Austrian Science Fund), \textit{Group Theory of Cyclic Cubic Fields}, MACOS 2024
\smallskip
\hrule

}

{\Large

\vspace{0.25in}

\noindent
I express my heartfelt gratitude to my

\bigskip
\begin{center}
\textbf{Coauthors:}
\end{center}

\bigskip
\noindent
\textbf{Siham Aouissi}, Universit\'e Moulay Ismail, Mekn\`es, Maroc,\\
\textbf{Bill Allombert}, Universit\'e de Bordeaux, France,\\
\textbf{Abderazak Soullami}, Universit\'e Sidi Mohammed Ben Abdellah, F\`es, Maroc.

\vspace{0.5in}

\noindent
My own ORCID: \texttt{https://orcid.org/my-orcid?orcid=0000-0002-4258-6455}.\\
SciProfiles: \texttt{https://sciprofiles.com/profile/3155986}.\\
ScopusPreview:\\ \texttt{https://www.scopus.com/authid/detail.uri?authorId=36696761400}.

\vspace{0.5in}


\begin{center}
\S\ 1. \textbf{Conductor, multiplicity, and 3-class rank}
\end{center}

\vspace{0.2in}

\noindent
A cyclic cubic number field \(K/\mathbb{Q}\)
whose \textbf{conductor} \(c=p_1\cdots p_t\)
is divisible by \(t\) primes \(p_i\equiv +1\,(\mathrm{mod}\,3)\)
or the prime power \(p_t=3^2=9\)
is member of a \textit{multiplet} \((K_1,\ldots,K_m)\)
with \(m=2^{t-1}\) non-isomorphic cyclic cubic fields
sharing the common \textit{discriminant} \(\mathrm{disc}(K_\mu)=c^2\),
for \(1\le\mu\le m\).
The number \(m=m(t)\) is called the \textbf{multiplicity} of the conductor \(c\),
but it only depends on the number \(t\).


\renewcommand{\arraystretch}{1.3}

\begin{table}[ht]
\caption{\textbf{Multiplicity} of conductors}
\label{tbl:Multiplicity}
\begin{center}

\begin{tabular}{|c||cccccc|}
\hline
 \(t\)       & \(0\) & \(1\) & \(2\) & \(3\) &  \(4\) &  \(5\) \\
\hline
 \(m= m(t)\) & \(0\) & \(2^0=1\) & \(2^1=2\) & \(2^2=4\) & \(2^3=8\) & \(2^4=16\) \\
 multiplet   & nilet &   singlet &   doublet &   quartet &     octet & hexadecuplet \\
\hline
\end{tabular}

\end{center}
\end{table}


\vspace{0.2in}

\noindent
\textit{Georges Gras} has shown that the \textbf{rank}
of the \(3\)-\textbf{class group}
\(\mathrm{Cl}_3(K)=\mathrm{Syl}_3\mathrm{Cl}(K)\)
of a cyclic cubic field \(K\) is bounded by
\(t-1\le\varrho\le 2\cdot(t-1)\),
in particular,\\
\(2 \le\varrho\le 4\) for \textit{quartets} \((K_1,\ldots,K_4)\), and
\(3\le\varrho\le 6\) for \textit{octets} \((K_1,\ldots,K_8)\).


\renewcommand{\arraystretch}{1.2}

\begin{table}[ht]
\caption{\textbf{Rank} of \(3\)-class groups}
\label{tbl:Rank}
\begin{center}

\begin{tabular}{|c||ccccc|}
\hline
 \(t\)     & \(1\) & \(2\) & \(3\) &  \(4\) &  \(5\) \\
\hline
 range     & \(0\le\varrho\le 0\) & \(1\le\varrho\le 2\) & \(2\le\varrho\le 4\) & \(3\le\varrho\le 6\) & \(4 \le\varrho\le 8\) \\
 multiplet &    singlet &   doublet &   quartet &     octet & hexadecuplet \\
\hline
\end{tabular}

\end{center}
\end{table}

}


{\normalsize

\vspace{0.25in}

\hrule
\smallskip
\noindent
Download of this presentation from \quad \texttt{http://www.algebra.at/DCM@MACOS2024Brasov.pdf}

}

\newpage

{\normalsize

\noindent
Daniel C. Mayer (Austrian Science Fund), \textit{Group Theory of Cyclic Cubic Fields}, MACOS 2024
\smallskip
\hrule

}

{\Large

\vspace{0.25in}

\begin{center}
\S\ 2. \textbf{Unramified abelian \(3\)-extensions}
\end{center}

\vspace{0.2in}

\noindent
The \textit{Hilbert \(3\)-class field tower} of a cyclic cubic field \(K\) is denoted by

\[
K < \mathrm{F}_3^1(K) < \mathrm{F}_3^2(K) < \ldots < \mathrm{F}_3^\infty(K), \text{ where }
\]

\bigskip
\noindent
\(\mathrm{F}_3^1(K)\) is the maximal \textbf{abelian} unramified \(3\)-extension of \(K\),\\
\(\mathrm{F}_3^2(K)\) is the maximal \textbf{metabelian} unramified \(3\)-extension of \(K\),\\ and so on, with increasing soluble length of the Galois group, until\\
\(\mathrm{F}_3^\infty(K)\) is the maximal unramified \textbf{pro-\(3\)-extension} of \(K\).

\vspace{0.25in}

\noindent
\textbf{Definition 1.}\\
By the \textbf{group theory} of \(K\) we understand
the identification of the automorphism group
\(M=\mathrm{Gal}(\mathrm{F}_3^2(K)/K)\)
of the \textit{maximal metabelian unramified \(3\)-extension} of \(K\),
which is called the \textbf{second \(3\)-class group} of \(K\).

\vspace{0.25in}

\noindent
A cyclic cubic field \(K\) with \(3\)-class rank \(\varrho\)
possesses \(U(\varrho)=\frac{3^\varrho-1}{3-1}\)
unramified cyclic cubic relative extensions \(E_u\), \(1\le u\le U(\varrho)\).

\bigskip
\noindent
\textbf{Definition 2.}\\
By the \textit{abelian type invariants} of \(K\)
we understand the family
\(\alpha(K)=\lbrack\mathrm{Cl}_3(E_u)\rbrack_{u=1}^{U(\varrho)}\),
and the \textit{capitulation type} of \(K\) is defined as the family
\(\varkappa(K)=(\ker(T_u))_{u=1}^{U(\varrho)}\),
where the \textit{transfer homomorphisms} from \(K\) to \(E_u\)
are denoted by \(T_u:\,\mathrm{Cl}_3(K)\to\mathrm{Cl}_3(E_u)\).

\vspace{0.25in}

\noindent
Each extension \(E_u\), \(1\le u\le U(\varrho)\),
with \(3\)-class rank \(\varrho_u\)
has \(U(\varrho_u)=\frac{3^{\varrho_u}-1}{3-1}\)
unramified cyclic cubic relative extensions \(E_{u,v}\), \(1\le v\le U(\varrho_u)\),
which are \textit{not necessarily abelian over} \(K\).

\bigskip
\noindent
\textbf{Definition 3.}\\
By the \textit{abelian type invariants of second order} of \(K\)
we understand the family
\[
\alpha_2(K)=(\mathrm{Cl}_3(E_u);\lbrack\mathrm{Cl}_3(E_{u,v})\rbrack_{v=1}^{U(\varrho_u)})_{u=1}^{U(\varrho)}.
\]

\vspace{0.25in}

\noindent
Since \textit{doublets} \((K_1,K_2)\) with \(t=2\), \(m=2\),
have been clarified completely with respect to \(M\)
in the doctoral thesis of \textit{Mohammed Ayadi},
we present brand new results for
\textit{quartets} with \(t=3\), \(m=4\),
and \textit{octets} with \(t=4\), \(m=8\).

We focus on \textit{octets} \((K_1,\ldots,K_8)\),
where \textit{all} components have elementary \textit{tri}cyclic \(3\)-class group
\(\mathrm{Cl}_3(K_\mu)\simeq (3,3,3)\), \(1\le\mu\le 8\),
and on \textit{quartets} \((K_1,\ldots,K_4)\),
where only a \textit{single} component has elementary \textit{tri}cyclic \(3\)-class group,
\(\mathrm{Cl}_3(K_1)\simeq (3,3,3)\),
and thus \(U(3)=\frac{3^3-1}{2}=13\) unramified cyclic cubic relative extensions
\((E_u)_{u=1}^{13}\).

}

{\normalsize

\vspace{0.25in}

\hrule
\smallskip
\noindent
Download of this presentation from \quad \texttt{http://www.algebra.at/DCM@MACOS2024Brasov.pdf}

}

\newpage

{\normalsize

\noindent
Daniel C. Mayer (Austrian Science Fund), \textit{Group Theory of Cyclic Cubic Fields}, MACOS 2024
\smallskip
\hrule

}

{\Large

\vspace{0.25in}

\begin{center}
\S\ 3. \textbf{Quartets with FOUR BIcyclic 3-class groups}
\end{center}

\vspace{0.2in}

\noindent
First we point out a most astonishing result, which concerns
\textit{quartets} \((K_1,\ldots,K_4)\) with \(m=4\), \(t=3\),
having \textit{four} elementary \textit{bi}cyclic \(3\)-class groups
\(\mathrm{Cl}_3(K_\mu)\simeq (3,3)\), \(1\le\mu\le 4\),
that is, a \textit{homogeneous} rank distribution \(\rho=(2^4)\).
Theorem 1 was derived in cooperation with my coauthor \textbf{Siham Aouissi} [2]:

\vspace{0.2in}

\noindent
\textbf{Theorem 1.}
\textit{If the prime divisors \(p,q,r\) of the conductor \(c=pqr\)
satisfy the following} \textbf{cubic residue conditions},
\textit{called} \textbf{graph III.8}
(abbreviated by the symbol
\(\lbrack p,q,r\rbrack_3=r\rightarrow p\leftrightarrow q\leftarrow r\),
and drawn in the following figure),
\[
\left(\frac{p}{q}\right)_3=\left(\frac{q}{p}\right)_3=1, \
\left(\frac{r}{p}\right)_3=1,\ \left(\frac{r}{q}\right)_3=1,\
\textit{ all other symbols }\ne 1,
\]
\textit{then all four \(3\)-class groups are elementary bicyclic
\(\mathrm{Cl}_3(K_\mu)\simeq (3,3)\), \(1\le\mu\le 4\),
and the second \(3\)-class group \(M\) of each component \(K_\mu\)
of the quartet \((K_1,\ldots,K_4)\) has} \textbf{coclass} \(\mathrm{cc}(M)\ge 2\),
\textit{that is, it} \textbf{cannot be of maximal nilpotency class}.

\begin{figure}[ht]
\label{fig:QuartetsIII8}

\setlength{\unitlength}{0.8cm}
\begin{picture}(6,4.5)(0,-5.5)

\put(2,-1.5){\makebox(0,0)[cc]{\(q\)}}
\put(2,-2){\circle*{0.2}}
\put(4,-1.5){\makebox(0,0)[cc]{\(p\)}}
\put(4,-2){\circle*{0.2}}

\put(2.2,-1.9){\vector(1,0){1.6}}
\put(3.8,-2.1){\vector(-1,0){1.6}}

\put(3,-4){\vector(-1,2){0.9}}
\put(3,-4){\vector(1,2){0.9}}
\put(4.8,-3){\makebox(0,0)[cc]{III.8}}

\put(3,-4){\circle*{0.2}}
\put(3,-4.5){\makebox(0,0)[cc]{\(r\)}}

\end{picture}

\end{figure}

\noindent
For \textit{all other situations} of cubic residue symbols,
\textit{maximal class}, \(\mathrm{cc}(M)=1\), \textit{is admissible} for the group \(M\),
and is confirmed by concrete numerical examples.

\vspace{0.2in}

However, the \textbf{statistical distribution} of the quartets
satisfying the assumptions of Theorem 1 is \textbf{very sparse}:
among the \(15\,851\) cyclic cubic fields with conductor \(c<10^5\),
there occur \(3132\), \(19,8\%\), in quartets.
Among these, \(2316\), \(73,9\%\), have four bicyclic \(3\)-class groups.
And among the latter, only \(28\), \(\mathbf{1.2\%}\), have \textbf{graph III.8},
with \textbf{minimal conductor} \(c=\mathbf{20\,293}=7\cdot 13\cdot 223\),
for which \(M\simeq\mathrm{SmallGroup}(729,i)\) with \(i\in\lbrace 37,38,39\rbrace\) [3],
and indeed \(\mathrm{cc}(M)=2\).

}

{\normalsize

\vspace{0.25in}

\hrule
\smallskip
\noindent
Download of this presentation from \quad \texttt{http://www.algebra.at/DCM@MACOS2024Brasov.pdf}

}

\newpage

{\normalsize

\noindent
Daniel C. Mayer (Austrian Science Fund), \textit{Group Theory of Cyclic Cubic Fields}, MACOS 2024
\smallskip
\hrule

}

{\Large

\vspace{0.25in}

\begin{center}
\S\ 4. \textbf{Quartets with ONE TRIcyclic 3-class group}
\end{center}

\vspace{0.15in}

\noindent
The next theorem concerns
\textit{quartets} \((K_1,\ldots,K_4)\) with \(m=4\), \(t=3\),
having \textit{three} elementary \textit{bi}cyclic \(3\)-class groups
\(\mathrm{Cl}_3(K_\mu)\simeq (3,3)\), \(2\le\mu\le 4\),
and a \textit{single} elementary \textit{tri}cyclic \(3\)-class group
\(\mathrm{Cl}_3(K_1)\simeq (3,3,3)\),
thus a \textit{heterogeneous} rank distribution \(\rho=(3,2^3)\).
It was proved together with coauthor \textbf{Bill Allombert} [1]:

\vspace{0.15in}

\noindent
\textbf{Theorem 2.}
\textit{If the prime divisors \(p,q,r\) of the conductor \(c=pqr\)
satisfy the following} \textbf{cubic residue conditions},
\textit{called} \textbf{graph I.2}
(abbreviated by the symbol
\(\lbrack p,q,r\rbrack_3=q\leftarrow p\rightarrow r\),
and drawn in the following figure),
\[
\left(\frac{p}{q}\right)_3=1,\ \left(\frac{p}{r}\right)_3=1, \
\textit{ all other symbols }\ne 1,
\]
\textit{then precisely one of the four \(3\)-class groups is elementary tricyclic
\(\mathrm{Cl}_3(K_1)\simeq (3,3,3)\),
and if the} \textbf{abelian type invariants (ATI)} \textit{of \(K_1\) are}
\[
\alpha(K_1)=\lbrack (9,3,3)^1,(9,9)^{12}\rbrack
\]
\textbf{of first order},
\textit{with rank distribution \(\rho(K_1)=(3^1,2^{12})\), and}
\begin{equation*}
\begin{aligned}
\alpha_2(K_1)=(
& \lbrack (9,3,3);(9,9,3),(9,3,3)^{12}\rbrack^1, \\
& \lbrack (9,9);(27,3,3)^3,(9,9,3)^1\rbrack^1, \\
& \lbrack (9,9);(9,9,3)^1,(9,3,3)^3\rbrack^2, \\
& \lbrack (9,9);(27,3,3)^1,(9,3,3)^3\rbrack^9)
\end{aligned}
\end{equation*}
\textbf{of second order},
\textit{then the second \(3\)-class group
\(M=\mathrm{Gal}(\mathrm{F}_3^2(K_1)/K_1)\)
of the component \(K_1\) with elementary tricyclic
\(3\)-class group \(\mathrm{Cl}_3(K_1)\simeq (3,3,3)\)
is the uniquely determined} \textbf{closed Andozhskii-Tsvetkov-group}
\(M\simeq\mathrm{SmallGroup}(6561,217710)\) [7]
\textit{with} \textbf{harmonically balanced capitulation (HBC)}, \(\varkappa\in S_{13}\),
\textit{a permutation}.

\begin{figure}[ht]
\label{fig:QuartetsI2}

\setlength{\unitlength}{0.8cm}
\begin{picture}(6,4.5)(0,-5.5)

\put(2,-1.5){\makebox(0,0)[cc]{\(r\)}}
\put(2,-2){\circle*{0.2}}
\put(4,-1.5){\makebox(0,0)[cc]{\(q\)}}
\put(4,-2){\circle*{0.2}}

\put(3,-4){\vector(-1,2){0.9}}
\put(3,-4){\vector(1,2){0.9}}
\put(4.8,-3){\makebox(0,0)[cc]{I.2}}

\put(3,-4){\circle*{0.2}}
\put(3,-4.5){\makebox(0,0)[cc]{\(p\)}}

\end{picture}

\end{figure}

\noindent
The \textbf{minimal conductor} is \(c=\mathbf{753\,787}=19\cdot 97\cdot 409\),
for the variant with rank distribution \(\rho(K_1)=(3^1,2^{12})\).
However, there exist two other variants,
and \(c=\mathbf{689\,347}=31\cdot 37\cdot 601\) is the minimal conductor
for the variant with rank distribution \(\rho(K_1)=(3^4,2^9)\)
and closed Andozhskii-Tsvetkov-group
\(M\simeq\mathrm{SmallGroup}(6561,i)\),
where \(i\in\lbrace 217713,217717\rbrace\) is not unique.
\textbf{Closed} means the relation rank
\(d_2(M)=\dim_{\mathbb{F}_3}H^2(M,\mathbb{F}_3)\)
coincides with the generator rank
\(d_1(M)=\dim_{\mathbb{F}_3}H^1(M,\mathbb{F}_3)=3\).

}

{\normalsize

\vspace{0.25in}

\hrule
\smallskip
\noindent
Download of this presentation from \quad \texttt{http://www.algebra.at/DCM@MACOS2024Brasov.pdf}

}

\newpage

{\normalsize

\noindent
Daniel C. Mayer (Austrian Science Fund), \textit{Group Theory of Cyclic Cubic Fields}, MACOS 2024
\smallskip
\hrule

}

{\Large

\vspace{0.25in}

\begin{center}
\S\ 5. \textbf{Octets with EIGHT TRIcyclic 3-class groups}
\end{center}

\vspace{0.2in}

\noindent
Now we come to the main result concerning
\textit{octets} \((K_1,\ldots,K_8)\) with \(m=8\), \(t=4\),
having \textit{eight} elementary \textit{tri}cyclic \(3\)-class groups
\(\mathrm{Cl}_3(K_\mu)\simeq (3,3,3)\), \(1\le\mu\le 8\),
\(\rho=(3^8)\).
It was established together with my coauthor \textbf{Abderazak Soullami}:

\vspace{0.2in}

\noindent
\textbf{Theorem 3.}
\textit{If the prime divisors \(p,q,r,s\) of the conductor \(c=pqrs\)
satisfy either of the} \textbf{cubic residue conditions},
\textit{called} \textbf{graph (v)}, \textit{respectively}, \textbf{graph (vii)}
(abbreviated by the symbol
\(\lbrack p,q,r,s\rbrack_3=p\rightarrow q\rightarrow r\rightarrow p\),
respectively
\(\lbrack p,q,r,s\rbrack_3=p\rightarrow q\rightarrow r\rightarrow s\rightarrow p\),
and drawn in the following figure),
\[
\left(\frac{p}{q}\right)_3=1,\ \left(\frac{q}{r}\right)_3=1, \
\left(\frac{r}{p}\right)_3=1,\
\textit{ all other symbols }\ne 1,
\]
\textit{respectively},
\[
\left(\frac{p}{q}\right)_3=1,\ \left(\frac{q}{r}\right)_3=1, \
\left(\frac{r}{s}\right)_3=1,\ \left(\frac{s}{p}\right)_3=1,\
\textit{ all other symbols }\ne 1,
\]
\textit{then all eight \(3\)-class groups are elementary tricyclic
\(\mathrm{Cl}_3(K_\mu)\simeq (3,3,3)\), \(1\le\mu\le 8\),
and if the} \textbf{abelian type invariants} \textit{are of the shape}
\(\alpha(K_\mu)=\lbrack\alpha,(3,3,3)^{12}\rbrack\),
\textit{then the second \(3\)-class group of \(K_\mu\), \(1\le\mu\le 8\), is}
\(M=\mathrm{Gal}(\mathrm{F}_3^2(K_\mu)/K_\mu)\simeq\)
\[
\begin{cases}
\mathrm{SmallGroup}(729,290) & \text{ if } \alpha=(3,3,3,3,3), \\
\mathrm{SmallGroup}(729,i),\ i\in\lbrace 291,293,294\rbrace & \text{ if } \alpha=(9,3,3,3), \\
\mathrm{SmallGroup}(2187,j),\ j\in\lbrace 2053,2055,2056\rbrace & \text{ if } \alpha=(9,3,3,3,3), \\
\mathrm{SmallGroup}(2187,k),\ k\in\lbrace 2058,\ldots,2065\rbrace & \text{ if } \alpha=(9,9,3,3)\qquad [3]
\end{cases}
\]
\textit{in dependence on the} \textbf{polarization} \(\alpha\).
\textit{The stabilization \((3,3,3)^{12}\) remains fixed}.

\vspace{0.2in}

\begin{figure}[ht]
\label{fig:TwoGraphs}

\setlength{\unitlength}{0.8cm}
\begin{picture}(16,4)(-14,-11.5)

\put(-9,-7){\makebox(0,0)[cc]{(v)}}
\put(-5,-7){\makebox(0,0)[cc]{(vii)}}
\multiput(-10,-7.5)(4,0){2}{\makebox(0,0)[cc]{\(s\)}}
\multiput(-8,-7.5)(4,0){2}{\makebox(0,0)[cc]{\(r\)}}
\multiput(-10,-8)(2,0){4}{\circle*{0.2}}
\put(-8,-8){\vector(-1,-1){1.8}}
\put(-4,-8){\vector(-1,0){1.8}}
\put(-6,-8){\vector(0,-1){1.8}}
\multiput(-8,-10)(4,0){2}{\vector(0,1){1.8}}

\multiput(-10,-10)(2,0){4}{\circle*{0.2}}
\multiput(-10,-10)(4,0){2}{\vector(1,0){1.8}}
\multiput(-10,-10.5)(4,0){2}{\makebox(0,0)[cc]{\(p\)}}
\multiput(-8,-10.5)(4,0){2}{\makebox(0,0)[cc]{\(q\)}}
\put(-9,-11){\makebox(0,0)[cc]{\(3\)-cycle}}
\put(-5,-11){\makebox(0,0)[cc]{\(4\)-cycle}}

\end{picture}

\end{figure}

\noindent
\textbf{Prototypes} with \textbf{minimal conductors} are
\(\mathbf{128\,583}=3^2\cdot 7\cdot 13\cdot 157\),
for the \(3\)-cycle graph (v), 
where \((9,3,3,3,3);(3,3,3,3,3)^2,(9,3,3,3)^5\)
is the distribution of polarizations among the octet,
and
\(\mathbf{188\,461}=7\cdot 13\cdot 19\cdot 109\),
for the \(4\)-cycle graph (vii),
with distribution of polarizations
\((9,9,3,3);(9,3,3,3)^7\).

}

{\normalsize

\vspace{0.25in}

\hrule
\smallskip
\noindent
Download of this presentation from \quad \texttt{http://www.algebra.at/DCM@MACOS2024Brasov.pdf}

}

\newpage

{\normalsize

\noindent
Daniel C. Mayer (Austrian Science Fund), \textit{Group Theory of Cyclic Cubic Fields}, MACOS 2024
\smallskip
\hrule

}

\vspace{0.25in}

{\normalsize

\begin{center}
\S\ 6. \textbf{System matrix and 3-class rank}
\end{center}

\vspace{0.25in}

\noindent
Let the conductor of the cyclic cubic field \(K\) be
\(c=p_1\cdots p_t\) with \(t\) prime divisors, 
either \(p_i\equiv +1\,(\mathrm{mod}\,3)\) when \(9\nmid c\)
or \(p_t=3^2\) when \(9\mid c\),
such that
the \textit{cubic residue symbols} are
\(\left(\frac{p_i}{p_j}\right)_3=\zeta_3^{a_{ij}}\)
for \(1\le i,j\le t\), \(i\ne j\),
where \(a_{ij}\in\lbrace -1,0,+1\rbrace\).

Let the \textit{Kummer generator} of \(K(\zeta_3)=\mathbb{Q}(\zeta_3,\alpha)\)
be
\(\alpha=c\cdot (a+3b\sqrt{-3})/2\in\mathbb{Q}(\zeta_3)\),
such that \(a,b\) satisfy the diophantine norm equation
\(4c=a^2+27b^2\) with \(a\in\mathbb{Z}\), \(b\in\mathbb{N}\), 
\(a\equiv +1\,(\mathrm{mod}\,3)\) when \(9\nmid c\),
but \(a=3a_0\), \(a_0\equiv +1\,(\mathrm{mod}\,3)\), \(b\equiv\pm 1\,(\mathrm{mod}\,3)\)
when \(9\mid c\),
and put \(v_i\equiv v_{\mathfrak{p}_{i1}}(\alpha)\equiv v_{\mathfrak{p}_{i2}}(\alpha)\,(\mathrm{mod}\,3)\) when \(p_i\equiv +1\,(\mathrm{mod}\,3)\),
in terms of the prime ideals
\(p_i\mathbb{Z}\lbrack\zeta_3\rbrack=\mathfrak{p}_{i1}\mathfrak{p}_{i2}\)
over the primes \(p_i\) which split completely in \(\mathbb{Q}(\zeta_3)\),
but \(v_t\equiv\frac{\alpha-1}{1-\zeta_3}\,(\mathrm{mod}\,\mathfrak{p}_t)\)
with \(3\mathbb{Z}\lbrack\zeta_3\rbrack=\mathfrak{p}_t^2\), i.e. \(\mathfrak{p}_t=(1-\zeta_3)\),
when \(p_t=3^2\).

In dependence on the number \(t\) of prime divisors of the conductor \(c\) of \(K\),
the \textbf{system matrix} \(M\)
is given in terms of the \textbf{cubic residue exponents}
\(a_{ij}\) for \(1\le i,j\le t\), \(i\ne j\),
and of the \textbf{valuations of the Kummer generator}
\(v_i\) for \(1\le i\le t\), by the following entries for \(1\le t\le 4\).

\begin{itemize}
\item
For \(t=1\):
\begin{equation}
\label{eqn:Matrix1}
M=
\begin{pmatrix}
0 \\
\end{pmatrix}
\end{equation}
\item
For \(t=2\):
\begin{equation}
\label{eqn:Matrix2}
M=
\begin{pmatrix}
-a_{12}v_2 &  a_{21}v_1 \\
 a_{12}v_2 & -a_{21}v_1 \\
\end{pmatrix}
\end{equation}
\item
\textbf{For} \(\mathbf{t=3}\):
\begin{equation}
\label{eqn:Matrix3}
M=
\begin{pmatrix}
-a_{12}v_2-a_{13}v_3 & a_{21}v_1 & a_{31}v_1 \\
a_{12}v_2 & -a_{21}v_1-a_{23}v_3 & a_{32}v_2 \\
a_{13}v_3 & a_{23}v_3 & -a_{31}v_1-a_{32}v_2 \\
\end{pmatrix}
\end{equation}
\item
\textbf{For} \(\mathbf{t=4}\):
\begin{equation}
\label{eqn:Matrix4}
M=
\begin{pmatrix}
-a_{12}v_2-a_{13}v_3-a_{14}v_4 & a_{21}v_1 & a_{31}v_1 & a_{41}v_1 \\
a_{12}v_2 & -a_{21}v_1-a_{23}v_3-a_{24}v_4 & a_{32}v_2 & a_{42}v_2 \\
a_{13}v_3 & a_{23}v_3 & -a_{31}v_1-a_{32}v_2-a_{34}v_4 & a_{43}v_3 \\
a_{14}v_4 & a_{24}v_4 & a_{34}v_4 & -a_{41}v_1-a_{42}v_2-a_{43}v_3 \\
\end{pmatrix}
\end{equation}
\end{itemize}

}

{\normalsize

\vspace{0.25in}

\hrule
\smallskip
\noindent
Download of this presentation from \quad \texttt{http://www.algebra.at/DCM@MACOS2024Brasov.pdf}

}

\newpage

{\normalsize

\noindent
Daniel C. Mayer (Austrian Science Fund), \textit{Group Theory of Cyclic Cubic Fields}, MACOS 2024
\smallskip
\hrule

}

{\Large

\vspace{0.15in}

\begin{center}
\S\ 7. \textbf{Tame octets}
\end{center}

\bigskip
\noindent
Now we define a kind of octet
all of whose sub-quartets belong to
Category III, Graphs \(1,\ldots,4\),
in the sense of [2].

\vspace{0.1in}

\noindent
\textbf{Definition 4.}\\
The graph of the conductor \(c=pqrs\)
of a \textbf{tame octet} \((K_1,\ldots,K_8)\)
must satisfy the following four \textit{axioms},
as illustrated in Figure
\ref{fig:SevenGraphs}
(\(u=\#\) of uni-directional edges).
\begin{enumerate}
\item
There are no bi-directional edges, \(b=0\).
\item
There are no attractive vertices, \(A=0\).
\item
There are no repulsive vertices, \(R=0\).
\item
For all sub graphs with three vertices
but without edges,\\
the \(\delta\)-invariant is not congruent to zero modulo \(3\),
\(\delta\not\equiv 0\,(\mathrm{mod}\,3)\).
\end{enumerate}

\vspace{0.10in}

\noindent
\textbf{Theorem A.}
\textit{An octet \((K_1,\ldots,K_8)\)
whose conductor \(c=pqrs\) satisfies the conditions in Definition \(4\)
has four sub-quartets with conductors
\(f\in\lbrace pqr,pqs,prs,qrs\rbrace\) belonging to
Category \(\mathrm{III}\), Graphs \(1,\ldots,4\).
Its \(4\)-graph is one of the seven graphs in Figure
\ref{fig:SevenGraphs}
(\(I=\#\) of isolated vertices).
}

\begin{proof}
The axioms discourage all graphs
different from Category III, Graphs \(1,\ldots,4\). \\
Axiom (1) eliminates sub graphs
in Category III, Graphs \(5,\ldots,9\). \\
Axiom (2) disables Graphs \(1\) and \(2\) of Category II. \\
Axiom (3) excludes Graph \(2\) of Category I. \\
Axiom (4) prohibits Graph \(1\) of Category I.
\end{proof}

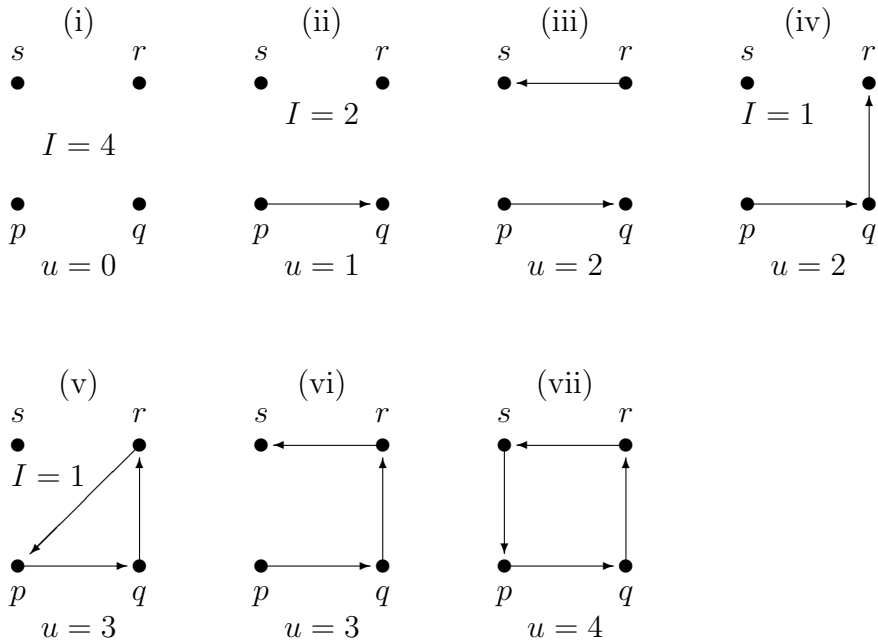
\begin{figure}[ht]
\caption{Seven graphs of tame octets}
\label{fig:SevenGraphs}

\setlength{\unitlength}{0.8cm}
\begin{picture}(16,10)(-10,-10.5)

\put(-9,-1){\makebox(0,0)[cc]{(i)}}
\put(-5,-1){\makebox(0,0)[cc]{(ii)}}
\put(-1,-1){\makebox(0,0)[cc]{(iii)}}
\put(3,-1){\makebox(0,0)[cc]{(iv)}}
\multiput(-10,-1.5)(4,0){4}{\makebox(0,0)[cc]{\(s\)}}
\multiput(-8,-1.5)(4,0){4}{\makebox(0,0)[cc]{\(r\)}}
\multiput(-10,-2)(2,0){8}{\circle*{0.2}}
\put(0,-2){\vector(-1,0){1.8}}
\put(4,-4){\vector(0,1){1.8}}

\put(-9,-3){\makebox(0,0)[cc]{\(I=4\)}}
\put(-5,-2.5){\makebox(0,0)[cc]{\(I=2\)}}
\put(2.5,-2.5){\makebox(0,0)[cc]{\(I=1\)}}

\multiput(-10,-4)(2,0){8}{\circle*{0.2}}
\multiput(-6,-4)(4,0){3}{\vector(1,0){1.8}}
\multiput(-10,-4.5)(4,0){4}{\makebox(0,0)[cc]{\(p\)}}
\multiput(-8,-4.5)(4,0){4}{\makebox(0,0)[cc]{\(q\)}}
\put(-9,-5){\makebox(0,0)[cc]{\(u=0\)}}
\put(-5,-5){\makebox(0,0)[cc]{\(u=1\)}}
\put(-1,-5){\makebox(0,0)[cc]{\(u=2\)}}
\put(3,-5){\makebox(0,0)[cc]{\(u=2\)}}

\put(-9,-7){\makebox(0,0)[cc]{(v)}}
\put(-5,-7){\makebox(0,0)[cc]{(vi)}}
\put(-1,-7){\makebox(0,0)[cc]{(vii)}}
\multiput(-10,-7.5)(4,0){3}{\makebox(0,0)[cc]{\(s\)}}
\multiput(-8,-7.5)(4,0){3}{\makebox(0,0)[cc]{\(r\)}}
\multiput(-10,-8)(2,0){6}{\circle*{0.2}}
\put(-8,-8){\vector(-1,-1){1.8}}
\put(-4,-8){\vector(-1,0){1.8}}
\put(-2,-8){\vector(0,-1){1.8}}
\put(0,-8){\vector(-1,0){1.8}}
\multiput(-8,-10)(4,0){3}{\vector(0,1){1.8}}

\put(-9.5,-8.5){\makebox(0,0)[cc]{\(I=1\)}}

\multiput(-10,-10)(2,0){6}{\circle*{0.2}}
\multiput(-10,-10)(4,0){3}{\vector(1,0){1.8}}
\multiput(-10,-10.5)(4,0){3}{\makebox(0,0)[cc]{\(p\)}}
\multiput(-8,-10.5)(4,0){3}{\makebox(0,0)[cc]{\(q\)}}
\put(-9,-11){\makebox(0,0)[cc]{\(u=3\)}}
\put(-5,-11){\makebox(0,0)[cc]{\(u=3\)}}
\put(-1,-11){\makebox(0,0)[cc]{\(u=4\)}}

\end{picture}

\end{figure}

}

{\normalsize

\vspace{0.25in}

\hrule
\smallskip
\noindent
Download of this presentation from \quad \texttt{http://www.algebra.at/DCM@MACOS2024Brasov.pdf}

}

\newpage

{\normalsize

\noindent
Daniel C. Mayer (Austrian Science Fund), \textit{Group Theory of Cyclic Cubic Fields}, MACOS 2024
\smallskip
\hrule

}

{\large

\vspace{0.25in}

\noindent
\textbf{Lemma.}\\
The distribution of the sub graphs
of the four sub quartets
of a tame octet is as follows.
\begin{enumerate}
\item
Graph (i): \\
Graph III.1 for \(pqr\),
Graph III.1 for \(pqs\),
Graph III.1 for \(prs\),
Graph III.1 for \(qrs\).
\item
Graph (ii): \\
Graph III.2 for \(pqr\),
Graph III.2 for \(pqs\),
Graph III.1 for \(prs\),
Graph III.1 for \(qrs\).
\item
Graph (iii): \\
Graph III.2 for \(pqr\),
Graph III.2 for \(pqs\),
Graph III.2 for \(prs\),
Graph III.2 for \(qrs\).
\item
Graph (iv): \\
Graph III.3 for \(pqr\),
Graph III.2 for \(pqs\),
Graph III.1 for \(prs\),
Graph III.2 for \(qrs\).
\item
Graph (v): \\
Graph III.4 for \(pqr\),
Graph III.2 for \(pqs\),
Graph III.2 for \(prs\),
Graph III.2 for \(qrs\).
\item
Graph (vi): \\
Graph III.3 for \(pqr\),
Graph III.2 for \(pqs\),
Graph III.2 for \(prs\),
Graph III.3 for \(qrs\).
\item
Graph (vii): \\
Graph III.3 for \(pqr\),
Graph III.3 for \(pqs\),
Graph III.3 for \(prs\),
Graph III.3 for \(qrs\).
\end{enumerate}

\noindent
\textbf{Proof.}\\
This can be seen by considering the four possible sub graphs with three vertices.

\vspace{0.25in}

\begin{center}
\S\ 8. \textbf{Minimal and maximal t-graphs of t prime divisors}
\end{center}

\bigskip
\noindent
Let \(\mathcal{G}_0\) be the graph
with all prime divisors of the conductor \(c\) as vertices
\(V=\lbrace p_1,\ldots,p_t\rbrace\),
and potential directed edges (arrows)
\(E_0=\lbrace p_i\to p_j \mid 1\le i,j\le t,\ i\ne j\rbrace\)
describing the cubic residue symbols
\(\left(\frac{p_i}{p_j}\right)_3=\zeta_3^{a_{ij}}\)
by the rule
\(p_i\to p_j\) iff \(\left(\frac{p_i}{p_j}\right)_3=1\) iff \(a_{ij}=0\).
For \textit{mutual cubic residues}, put briefly
\(p_i\leftrightarrow p_j\) iff \(a_{ij}=a_{ji}=0\).
For all graphs \(\mathcal{G}=(V,E)\) in the sequel,
the set \(V\) of vertices remains fixed,
only the sets \(E\) of directed edges varies.

\bigskip
\noindent
\textbf{Definition 5.}\\
Let \(\mathcal{G}_m\) be a \textit{minimal} graph,
and \(\mathcal{G}_M\) be a \textit{maximal} graph
such that \(\mathcal{G}_m<\mathcal{G}_M\),
in the sense that
\(E_m\subset E_M\).
An arrow \(e\in E_0\) is called
\begin{equation}
\label{eqn:Graphs}
\begin{cases}
mandatory, \\
admissible, \\
forbidden,
\end{cases}
\quad \text{ iff } \quad
\begin{cases}
e\in E_m, \\
e\in E_M\setminus E_m, \\
e\in E_0\setminus E_M.
\end{cases}
\end{equation}
Then, a graph with the property 
\(\mathcal{G}_m\le\mathcal{G}\le\mathcal{G}_M\)
contains all mandatory arrows,
may contain admissible arrows,
but does certainly not contain any forbidden arrows.

}

{\normalsize

\vspace{0.25in}

\hrule
\smallskip
\noindent
Download of this presentation from \quad \texttt{http://www.algebra.at/DCM@MACOS2024Brasov.pdf}

}

\newpage

{\normalsize

\noindent
Daniel C. Mayer (Austrian Science Fund), \textit{Group Theory of Cyclic Cubic Fields}, MACOS 2024
\smallskip
\hrule

}

{\Large

\vspace{0.25in}

\begin{center}
\S\ 9. \textbf{Rank distribution of 4-graphs with four prime divisors}
\end{center}

\vspace{0.25in}

\noindent
Motivated by the proof of the Theorem C for \(t=3\) in the Appendix,
the next theorem shows that similar techniques can be applied
to \(t=4\), for which
\textbf{no statements are available in the literature up to now}. 

\bigskip
\noindent
\textbf{Theorem B.}
Let \(t=4\) and \(c=pqrs\).
Then the octet \(K_1,\ldots,K_8\)
possesses \textbf{homogeneous} rank distribution
with \textbf{elementary tricyclic} \(3\)-class groups
\(\mathrm{Cl}_3(K_\mu)=(3,3,3)\), \(1\le\mu\le 8\),
when the graph \(\mathcal{G}\) of \(p,q,r,s\)
belongs to the range
\(\mathcal{G}_m\le\mathcal{G}\le\mathcal{G}_M\), where

\begin{figure}[ht]
\label{fig:Octets}

\setlength{\unitlength}{0.8cm}
\begin{picture}(16,5)(-10,-5.5)

\put(-9,-1){\makebox(0,0)[cc]{\(\mathcal{G}_m\)}}
\put(-3,-1){\makebox(0,0)[cc]{\(\mathcal{G}\)}}
\put(1,-1){\makebox(0,0)[cc]{\(\cdots\)}}
\put(3,-1){\makebox(0,0)[cc]{\(\mathcal{G}_M\)}}
\multiput(-10,-1.5)(4,0){4}{\makebox(0,0)[cc]{\(s\)}}
\multiput(-8,-1.5)(4,0){4}{\makebox(0,0)[cc]{\(r\)}}
\multiput(-10,-2)(2,0){8}{\circle*{0.2}}

\multiput(-10,-4)(4,0){3}{\vector(0,1){1.8}}
\multiput(-10,-4)(4,0){3}{\vector(1,1){1.8}}
\multiput(-8,-4)(4,0){3}{\vector(-1,1){1.8}}

\multiput(-6,-2)(4,0){3}{\vector(1,0){1.8}}

\put(4,-4){\vector(-1,0){1.8}}
\multiput(0,-2)(4,0){2}{\vector(0,-1){1.8}}

\put(1.9,-3.8){\vector(0,1){1.6}}
\put(2.1,-2.2){\vector(0,-1){1.6}}
\put(2.2,-3.8){\vector(1,1){1.4}}
\put(3.9,-2.4){\vector(-1,-1){1.4}}
\put(3.5,-3.8){\vector(-1,1){1.4}}
\put(2.4,-2.4){\vector(1,-1){1.4}}

\multiput(-10,-4)(2,0){8}{\circle*{0.2}}
\multiput(-10,-4.5)(4,0){4}{\makebox(0,0)[cc]{\(p\)}}
\multiput(-8,-4.5)(4,0){4}{\makebox(0,0)[cc]{\(q\)}}
\put(-9,-5){\makebox(0,0)[cc]{\(u=3\)}}
\put(-5,-5){\makebox(0,0)[cc]{\(u=4\)}}
\put(-1,-5){\makebox(0,0)[cc]{\(u=5\)}}
\put(3,-5){\makebox(0,0)[cc]{\(u=3,b=3\)}}

\end{picture}

\end{figure}

\begin{proof}
Since the sum of all rows of the \((4\times 4)\)-system matrix \(M\) in Formula
(4)
is zero, the first row can be cancelled.
It is the negative sum of the remaining rows.
The resulting \((3\times 4)\)-matrix
\[
\begin{pmatrix}
a_{12}v_2 & -(a_{21}v_1+a_{23}v_3+a_{24}v_4) & a_{32}v_2 & a_{42}v_2 \\
a_{13}v_3 & a_{23}v_3 & -(a_{31}v_1+a_{32}v_2+a_{34}v_4) & a_{43}v_3 \\
a_{14}v_4 & a_{24}v_4 & a_{34}v_4 & -(a_{41}v_1+a_{42}v_2+a_{43}v_3) \\
\end{pmatrix}
\]
is reduced to echelon form
by the cubic residue conditions in the graph \(\mathcal{G}\),
namely the mandatory arrows \(a_{13}=a_{14}=a_{24}=0\):
\[
\begin{pmatrix}
a_{12}v_2 & -(a_{21}v_1+a_{23}v_3+a_{24}v_4) & a_{32}v_2 & a_{42}v_2 \\
0         & a_{23}v_3 & -(a_{31}v_1+a_{32}v_2+a_{34}v_4) & a_{43}v_3 \\
0         & 0         & a_{34}v_4 & -(a_{41}v_1+a_{42}v_2+a_{43}v_3) \\
\end{pmatrix}
\]
The forbidden arrows \(a_{12},a_{23},a_{34}\in\lbrace -1,+1\rbrace\)
generate the pivot elements in the echelon form.
Thus, \(r=t-1=3\) and \(\varrho=2(t-1)-r=6-3=3\),
according to section \S\ 1.
Finally, it is well known
[2]
that
\(\varrho=3\) for \(t=4\) implies elementary
\(\mathrm{Cl}_3(K_\mu)=(3,3,3)\), for all \(1\le\mu\le 8\).
\end{proof}

}

{\normalsize

\vspace{0.25in}

\hrule
\smallskip
\noindent
Download of this presentation from \quad \texttt{http://www.algebra.at/DCM@MACOS2024Brasov.pdf}

}

\newpage

{\normalsize

\noindent
Daniel C. Mayer (Austrian Science Fund), \textit{Group Theory of Cyclic Cubic Fields}, MACOS 2024
\smallskip
\hrule

}

{\Large

\vspace{0.25in}

\noindent
The following experimental results
have been computed with MAGMA
[4,5,6],
using the class field theoretic routines by \textit{Claus Fieker}, 2001.

\bigskip
\noindent
\textbf{Example.}
\label{exm:Homogeneous}
The conditions of Theorem C
are satisfied by octets \(K_1,\ldots,K_8\) with \(t=4\),
and the following conductors of the form \(c=pqrs\):
\begin{equation}
\label{eqn:Rank3}
c\in\lbrace 104\,013, 144\,781, 158\,067, 181\,503\rbrace,
\end{equation}
where
\(\mathbf{104\,013}=3^2\cdot 7\cdot 13\cdot 127\),
\(\mathbf{144\,781}=7\cdot 13\cdot 37\cdot 43\),
\(\mathbf{158\,067}=3^2\cdot 7\cdot 19\cdot 127\),
\(\mathbf{181\,503}=3^2\cdot 7\cdot 43\cdot 67\),
and all \(4\)-graphs are isomorphic to
the minimal graph \(\mathcal{G}_m\).
Thus, they possess homogeneous rank distribution \(\rho(c)=(3^8)\),
and elementary tricyclic
\(\mathrm{Cl}_3(K_\mu)=(3,3,3)\), for all \(1\le\mu\le 8\).

\vspace{0.25in}

\begin{center}
\S\ 10. \textbf{Tame homogeneous octets}
\end{center}

\vspace{0.25in}

\noindent
In the following table,
further results of the present article are summarized.
By means of the class field theoretic routines by \textit{Fieker},
which are implemented in MAGMA,
we have computed the \(13\) unramified abelian extensions \(E_u/K_\mu\), \(1\le u\le 13\),
the logarithmic abelian type invariants \((\alpha_\mu)_{\mu=1}^8\) of their \(3\)-class groups \(\mathrm{Cl}_3(E_u)\)
and the capitulation type  \((\varkappa_\mu)_{\mu=1}^8\) of their transfer kernels \(\ker(T_u)\),
for each member \(K_\mu\) of tame homogeneous octets \((K_1,\ldots,K_8)\)
with conductors in the range \(c\le 200\,000\) without gaps,
and some isolated cases beyond.
This \textbf{Artin pattern} \((\alpha,\varkappa)\) is used as a search pattern
for a database query in the SmallGroups Library [3],
which identifies the second \(3\)-class group
\(M_\mu=\mathrm{Gal}(\mathrm{F}_3^2(K_\mu)/K_\mu)\).

}

\renewcommand{\arraystretch}{1.5}

\begin{table}[ht]
\caption{Invariants of tame homogeneous octets}
\label{tbl:Tame}
\begin{center}

{\tiny

\begin{tabular}{|r|r||c|c||c|c||c|}
\hline
 No    & \(c\)        & \(p,q,r,s\)     & Graph & \((\alpha_\mu)_{\mu=1}^8\)                        & \((\varkappa_\mu)_{\mu=1}^8\)     & \((M_\mu)_{\mu=1}^8\) \\
\hline
 \(1\) &  \(30\,303\) & \(9,7,13,37\)   & (iii) & \(2211,(111)^{12}\), \((2111,(111)^{12})^7\)  & \(LO^{12}\), \((PO^{12})^7\)  & \(\langle 2187,2058\rangle\), \(\langle 729,291\rangle^7\) \\
 \(2\) &  \(51\,471\) & \(9,7,19,43\)   & (vi)  & \((2111,(111)^{12})^8\)                       & \((PO^{12})^8\)               & \(\langle 729,291\rangle^8\) \\
 \(3\) &  \(80\,199\) & \(9,7,19,67\)   & (vi)  & \(2211,(111)^{12}\), \(11111,(111)^{12}\), \((2111,(111)^{12})^6\) &
 \(LO^{12}\), \((PO^{12})^7\)  & \(\langle 2187,2058\rangle\), \(\langle 729,290\rangle\), \(\langle 729,291\rangle^6\) \\
 \(4\) & \(128\,583\) & \(9,7,13,157\)  & (v)   & \(21111,(111)^{12}\), \((11111,(111)^{12})^2\), \((2111,(111)^{12})^5\) &
 \(LO^{12}\), \((PO^{12})^7\)  & \(\langle 2187,2053\rangle\), \(\langle 729,290\rangle^2\), \(\langle 729,291\rangle^5\) \\
 \(5\) & \(162\,981\) & \(9,7,13,199\)  & (iii) & \((2111,(111)^{12})^8\)                       & \((PO^{12})^8\)               & \(\langle 729,291\rangle^8\) \\
 \(6\) & \(172\,081\) & \(7,13,31,61\)  & (vi)  & \(21111,(111)^{12}\), \(11111,(111)^{12}\), \((2111,(111)^{12})^6\) &
 \(LO^{12}\), \((PO^{12})^7\)  & \(\langle 2187,2053\rangle\), \(\langle 729,290\rangle\), \(\langle 729,291\rangle^6\) \\
 \(7\) & \(188\,461\) & \(7,13,19,109\) & (vii) & \(2211,(111)^{12}\), \((2111,(111)^{12})^7\)  & \(LO^{12}\), \((PO^{12})^7\)  & \(\langle 2187,2058\rangle\), \(\langle 729,291\rangle^7\) \\
 \(8\) & \(294\,903\) & \(9,7,31,151\)  & (iv)  & \(11111,(111)^{12}\), \((2111,(111)^{12})^7\) & \((PO^{12})^8\)               & \(\langle 729,290\rangle\), \(\langle 729,291\rangle^7\) \\
 \(9\) & \(397\,449\) & \(9,13,43,79\)  & (ii)  & \((2111,(111)^{12})^8\)                       & \((PO^{12})^8\)               & \(\langle 729,291\rangle^8\) \\
\(10\) & \(447\,237\) & \(9,7,31,229\)  & (i)   & \((2111,(211)^3,(111)^9)^8\)                  & \(((P_i)_{i=1}^4O^9)^8\)      & \(\langle 2187,814\rangle^8\) \\
\hline
\end{tabular}

}

\end{center}
\end{table}

{\normalsize

\vspace{0.25in}

\hrule
\smallskip
\noindent
Download of this presentation from \quad \texttt{http://www.algebra.at/DCM@MACOS2024Brasov.pdf}

}

\newpage

{\normalsize

\noindent
Daniel C. Mayer (Austrian Science Fund), \textit{Group Theory of Cyclic Cubic Fields}, MACOS 2024
\smallskip
\hrule

}

{\Large

\vspace{0.25in}

\begin{center}
Appendix. \textbf{Rank distribution of 3-graphs with three prime divisors}
\end{center}

\vspace{0.25in}

\noindent
The following theorem
proves a dominant part of 
Proposition VI.5, pp. 20--21, by Georges Gras, 1973,
which is given without proof,
justified by the explanation that
the proof is straight-forward but fastidious.
In the terminology of
[2],
the following theorem clarifies the
\textbf{Graphs 2,3,4,5,6,7,9 of Category III}.

\bigskip
\noindent
\textbf{Theorem C.}
\textit{Let \(t=3\) and \(c=pqr\).
Then the quartet \(K_1,\ldots,K_4\) 
possesses} \textbf{homogeneous} \textit{rank distribution
with} \textbf{elementary bicyclic} \textit{\(3\)-class groups
\(\mathrm{Cl}_3(K_i)=(3,3)\), \(1\le i\le 4\),
when the graph \(\mathcal{G}\) of \(p,q,r\)
belongs to the range
\(\mathcal{G}_m\le\mathcal{G}\le\mathcal{G}_M\), where}

\begin{figure}[ht]
\label{fig:Quartets}

\setlength{\unitlength}{0.8cm}
\begin{picture}(16,5)(-10,-5.5)

\put(-9,-1){\makebox(0,0)[cc]{\(\mathcal{G}_m\)}}
\put(-3,-1){\makebox(0,0)[cc]{\(\mathcal{G}\)}}
\put(1,-1){\makebox(0,0)[cc]{\(\cdots\)}}
\put(3,-1){\makebox(0,0)[cc]{\(\mathcal{G}_M\)}}
\multiput(-10,-1.5)(4,0){4}{\makebox(0,0)[cc]{\(p\)}}
\multiput(-8,-1.5)(4,0){4}{\makebox(0,0)[cc]{\(q\)}}
\multiput(-10,-2)(2,0){8}{\circle*{0.2}}

\multiput(-10,-2)(4,0){3}{\vector(1,-2){0.9}}
\put(-8,-3){\makebox(0,0)[cc]{III.2}}

\multiput(-4,-2)(4,0){3}{\vector(-1,0){1.8}}
\put(-4,-3){\makebox(0,0)[cc]{III.3}}

\multiput(-1,-4)(4,0){2}{\vector(1,2){0.9}}
\put(0.2,-3){\makebox(0,0)[cc]{III.4}}

\put(2.8,-3.8){\vector(-1,2){0.8}}
\put(2.2,-2.2){\vector(1,-2){0.8}}
\put(4.2,-3){\makebox(0,0)[cc]{III.9}}

\multiput(-9,-4)(4,0){4}{\circle*{0.2}}
\multiput(-9,-4.5)(4,0){4}{\makebox(0,0)[cc]{\(r\)}}
\put(-9,-5){\makebox(0,0)[cc]{\(u=1\)}}
\put(-5,-5){\makebox(0,0)[cc]{\(u=2\)}}
\put(-1,-5){\makebox(0,0)[cc]{\(u=3\)}}
\put(3,-5){\makebox(0,0)[cc]{\(u=2,b=1\)}}

\end{picture}

\end{figure}

\begin{proof}
Since the sum of all rows of the \((3\times 3)\)-system matrix \(M\) in Formula
(3)
is zero, the first row can be cancelled.
It is the negative sum of the remaining rows.
The resulting \((2\times 3)\)-matrix
\[
\begin{pmatrix}
a_{12}v_2 & -(a_{21}v_1+a_{23}v_3) & a_{32}v_2 \\
a_{13}v_3 & a_{23}v_3 & -(a_{31}v_1+a_{32}v_2) \\
\end{pmatrix}
\]
is reduced to echelon form
by the cubic residue conditions in the graph \(\mathcal{G}\),
namely the mandatory arrow \(a_{13}=0\):
\[
\begin{pmatrix}
a_{12}v_2 & -(a_{21}v_1+a_{23}v_3) & a_{32}v_2 \\
0         & a_{23}v_3 & -(a_{31}v_1+a_{32}v_2) \\
\end{pmatrix}
\]
The forbidden arrows \(a_{12},a_{23}\in\lbrace -1,+1\rbrace\)
generate the pivot elements in the echelon form.
Thus, \(r=t-1=2\) and \(\varrho=2(t-1)-r=4-2=2\),
according to section \S\ 1.
Finally, it is well known
[2]
that
\(\varrho=2\) for \(t=3\) implies elementary
\(\mathrm{Cl}_3(K_\mu)=(3,3)\), for all \(1\le\mu\le 4\).
\end{proof}

}

{\normalsize

\vspace{0.25in}

\hrule
\smallskip
\noindent
Download of this presentation from \quad \texttt{http://www.algebra.at/DCM@MACOS2024Brasov.pdf}

}

\newpage

{\normalsize

\noindent
Daniel C. Mayer (Austrian Science Fund), \textit{Group Theory of Cyclic Cubic Fields}, MACOS 2024
\smallskip
\hrule

}

{\Large

\vspace{0.25in}

\noindent
The following theorem
supplements another part of 
Proposition VI.5, pp. 20--21, by Georges Gras, 1973,
which is given without proof.
In the terminology of
[2],
the following theorem clarifies the
\textbf{Graphs 1, 2 of Category I},
the
\textbf{Graphs 1, 2 of Category II},
and the
\textbf{Graph 1 of Category III}.

\bigskip
\noindent
\textbf{Theorem D.}
Let \(t=3\), \(c=pqr\) and \(\delta=a_{12}a_{23}a_{31}-a_{13}a_{32}a_{21}\).
Then the quartet \(K_1,\ldots,K_4\) 
has
\begin{itemize}
\item
\textbf{inhomogeneous} rank distribution
with a single tricyclic \(3\)-class group
\(\mathrm{Cl}_3(K_1)\) of rank three,
and \textbf{three elementary bicyclic} \(3\)-class groups
\(\mathrm{Cl}_3(K_\mu)=(3,3)\), \(2\le\mu\le 4\),
when the graph \(\mathcal{G}\) of \(p,q,r\)
has no edges and \(\delta\equiv 0\,(\mathrm{mod}\,3)\)
(\textbf{Graph I.1});
\item
\textbf{homogeneous} rank distribution
with \textbf{four elementary bicyclic} \(3\)-class groups
\(\mathrm{Cl}_3(K_\mu)=(3,3)\), \(1\le\mu\le 4\),
when the graph \(\mathcal{G}\) of \(p,q,r\)
has no edges and \(\delta\not\equiv 0\,(\mathrm{mod}\,3)\)
(\textbf{Graph III.1});
\item
\textbf{inhomogeneous} rank distribution
with a single tricyclic \(3\)-class group
\(\mathrm{Cl}_3(K_1)\) of rank three,
and \textbf{three elementary bicyclic} \(3\)-class groups
\(\mathrm{Cl}_3(K_\mu)=(3,3)\), \(2\le\mu\le 4\),
when the graph \(\mathcal{G}\) of \(p,q,r\)
has precisely one universally repelling vertex
(\textbf{Graph I.2});
\item
\textbf{inhomogeneous} rank distribution
with two tricyclic \(3\)-class group
\(\mathrm{Cl}_3(K_\mu)\), \(1\le\mu\le 2\), of rank three,
and \textbf{two elementary bicyclic} \(3\)-class groups
\(\mathrm{Cl}_3(K_\mu)=(3,3)\), \(3\le\mu\le 4\),
when the graph \(\mathcal{G}\) of \(p,q,r\)
has either precisely one universally attracting vertex
(\textbf{Graph II.1})
or a universally repelling vertex additionally
(\textbf{Graph II.2}).
\end{itemize}

\begin{figure}[ht]
\label{fig:Quartets2}

\setlength{\unitlength}{0.8cm}
\begin{picture}(16,3.5)(-10,-5)

\multiput(-10,-1.5)(4,0){4}{\makebox(0,0)[cc]{\(p\)}}
\multiput(-8,-1.5)(4,0){4}{\makebox(0,0)[cc]{\(q\)}}
\multiput(-10,-2)(2,0){8}{\circle*{0.2}}

\put(-9.3,-3){\makebox(0,0)[cc]{III.1 and}}
\put(-8,-3){\makebox(0,0)[cc]{I.1}}

\multiput(-6,-2)(8,0){2}{\vector(1,-2){0.9}}
\multiput(-6,-2)(4,0){3}{\vector(1,0){1.8}}
\put(-4,-3){\makebox(0,0)[cc]{I.2}}

\put(-1,-4){\vector(1,2){0.9}}
\put(0.2,-3){\makebox(0,0)[cc]{II.1}}

\put(4,-2){\vector(-1,-2){0.9}}
\put(4.2,-3){\makebox(0,0)[cc]{II.2}}

\multiput(-9,-4)(4,0){4}{\circle*{0.2}}
\multiput(-9,-4.5)(4,0){4}{\makebox(0,0)[cc]{\(r\)}}
\put(-9,-5){\makebox(0,0)[cc]{\(u=0\)}}
\put(-5,-5){\makebox(0,0)[cc]{\(u=2\)}}
\put(-1,-5){\makebox(0,0)[cc]{\(u=2\)}}
\put(3,-5){\makebox(0,0)[cc]{\(u=3\)}}

\end{picture}

\end{figure}

}

{\normalsize

\vspace{0.25in}

\hrule
\smallskip
\noindent
Download of this presentation from \quad \texttt{http://www.algebra.at/DCM@MACOS2024Brasov.pdf}

}

\newpage

{\normalsize

\noindent
Daniel C. Mayer (Austrian Science Fund), \textit{Group Theory of Cyclic Cubic Fields}, MACOS 2024
\smallskip
\hrule

}

{\Large

\vspace{0.25in}

\begin{proof}
These results are proved computationally.
Nine nested loops are executed,
each over two values in \(\lbrace 1,2\rbrace\).
The first six loops over
\(a_{12},a_{13},a_{21},a_{23},a_{31},a_{32}\)
iterate the exponents
of cubic residue symbols
between the prime(power)s \(p,q,r\)
dividing the conductor \(c=p\cdot q\cdot r\)
of the cyclic cubic field \(K\).
The other three loops iterate the valuations
\(v_1,v_2,v_3\)
of the Kummer generator \(\alpha\) of \(K\).
For each of the \(2^9=512\) nonets,
the rank \(r\) of the system matrix \(M\) in Formula
(3)
is computed over the finite field \(\mathbb{F}_3\)
with MAGMA
[4,5,6].

For Graph 1 of Category I,
the constraint \(\delta\equiv 0\,(\mathrm{mod}\,3)\)
for the \(\delta\)-invariant
\(\delta=a_{12}a_{23}a_{31}-a_{13}a_{32}a_{21}\)
reduces the system to \(256\) nonets.
\(64\) of them yield \(r=1\), \(\varrho=2(t-1)-r=2\cdot 2-r=4-1=3\),
and \(192=3\cdot 64\) yield \(r=2\), \(\varrho=4-2=2\),
and elementary bicyclic \(3\)-class group.

For Graph 1 of Category III,
the constraint \(\delta\not\equiv 0\,(\mathrm{mod}\,3)\)
produces the remaining \(256\) nonets
with uniform rank \(r=2\), \(\varrho=4-2=2\),
and elementary bicyclic \(3\)-class group.

For Graph 2 of Category I,
we assume without loss of generality
that \(a_{12}=a_{13}=0\),
which narrows down to seven nested loops.
Among the \(2^7=128\) septets,
\(32\) yield \(r=1\), \(\varrho=4-1=3\),
and \(96=3\cdot 32\) yield \(r=2\), \(\varrho=4-2=2\),
and elementary bicyclic \(3\)-class group.

For Graph 1 of Category II,
we assume without loss of generality
that \(a_{12}=a_{32}=0\),
which reduces to seven nested loops.
Among the \(2^7=128\) septets,
\(64\) yield \(r=1\), \(\varrho=4-1=3\),
and \(64\) yield \(r=2\), \(\varrho=4-2=2\),
and elementary bicyclic \(3\)-class group.

For Graph 2 of Category II,
we assume without loss of generality
that \(a_{12}=a_{13}=a_{23}=0\),
which narrows down to six nested loops.
Among the \(2^6=64\) sextets,
\(32\) yield \(r=1\), \(\varrho=4-1=3\),
and \(32\) yield \(r=2\), \(\varrho=4-2=2\),
and elementary bicyclic \(3\)-class group.
\end{proof}

}

{\small

\vspace{0.2in}

\begin{center}
\textbf{References.}
\end{center}

\begin{enumerate}
\item[\lbrack 1\rbrack]
B. Allombert, D. C. Mayer,
\textit{Corps de nombres cubiques cycliques ayant une capitulation harmonieusement \'equilibr\'ee}, Publ. Math. Besan\c con.
\item[\lbrack 2\rbrack]
S. Aouissi, D. C. Mayer,
\textit{A group theoretic approach to cyclic cubic fields},
Mathematics
\textbf{12} (2024), nr. 126 , MDPI.
\item[\lbrack 3\rbrack]
H. U. Besche, B. Eick and E. A. O'Brien,
\textit{The SmallGroups Library --- a Library of Groups of Small Order},
2005,
an accepted and refereed GAP package, available also in MAGMA.
\item[\lbrack 4\rbrack]
W. Bosma, J. Cannon, and C. Playoust,
\textit{The Magma algebra system. I. The user language},
J. Symbolic Comput.
\textbf{24}
(1997),
235--265.
\item[\lbrack 5\rbrack]
W. Bosma, J. J. Cannon, C. Fieker, and A. Steel (eds.),
\textit{Handbook of Magma functions}
(Edition 2.28,
Sydney,
2024).
\item[\lbrack 6\rbrack]
MAGMA Developer Group,
MAGMA \textit{Computational Algebra System},
Version 2.28-9,
Sydney,
2024,
available from
\texttt{http://magma.maths.usyd.edu.au}.
\item[\lbrack 7\rbrack]
MAGMA Developer Group,
MAGMA \textit{Data for groups of order \(3^8\)},
\texttt{data3to8.tar.gz},
Sydney,
2012,
available from
\texttt{http://magma.maths.usyd.edu.au}.
\end{enumerate}

}

{\normalsize

\vspace{0.2in}

\hrule
\noindent
Download of this presentation from \quad \texttt{http://www.algebra.at/DCM@MACOS2024Brasov.pdf}

}


\end{document}